\documentclass[12pt]{amsart}
\usepackage{amsfonts}
\usepackage{amssymb}
\theoremstyle{definition}

\theoremstyle{remark}

\newcommand{\R}{\mathbb R}
\newcommand{\C}{\mathbb C}

\begin{document}

\centerline{\large\bf MATHEMATICS AND EDUCATION IN MATHEMATICS, 1986}
\centerline{\bf Proceedings of the Fifteenth Spring Conference of the Union of}
\centerline{\bf Bulgarian Mathematicians}
\centerline{\bf Sunny Beach, April 2-6, 1986}

\vspace{0.5in}
\centerline{\large ON THE BOUNDEDNESS OF THE SECTIONAL CURVATURE}
\centerline{\large OF ALMOST HERMITIAN MANIFOLDS }

\vspace{0.2in}
\centerline{\large Adrijan V. Borisov, Ognian T. Kassabov, }

\vspace{0.2in}
{\sl  In this paper we study the boundedness of the sectional curvatures
(holomorphic and antiholo\-morphic) of almost Hermitian manifolds of definite or
indefinite metrics. We prove that an almost Hermitian manifold \ $M$ \ of indefinite metric 
is of pointwise constant holomorphic sectional curvature  if the holomorphic sectional
curvature is bounded from above and from below. If the antiholomorphic sectional curvature 
is bounded either from above or from below, then \ $M$ \ is of pointwise constant 
antiholomorphic sectional curvature. Similar results are obtained for almost Hermitian 
manifolds of definite metric.}

\vspace{0.2in}
Examining indefinite Riemannian manifolds with bounded sectional curvature, R. S. Kulkarni [6]
proved that such a manifold is of constant sectional curvature. This theorem was specified by
M. Dajczer and K. Nomizu [3]. K\"ahler analogous of the theorem of R. S. Kulkarni was obtained 
in [1] (see also [5]). Now we show that the results in [1], [5] can be extended for an 
arbitrary almost Hermitian manifold. On the other hand, G. Ganchev and V. Mihova [4] used 
the complexification of the tangent space of a Riemannian manifold to examine the behaviour
of the sectional curvature. Using this approach we study the boundedness of the
sectional curvature in the complexification of the tangent space of an almost Hermitian
manifold of definite metric. 

\vspace{0.1in}
1. \underline{ Introduction}. Let \ $M$ \ be an almost Hermitian manifold with (definite or indefinite) 
metric tensor \ $g$ \ and almost complex structure \ $J$, i.e. \ $J^2=-id$, $g(X,Y)=g(JX,JY)$.
A plane \ $\alpha$ \ in \ $T_pM$ \ or in the complexification \ $(T_pM)^\C=T_pM\otimes_\R \C$ \ is
said to be holomorphic (resp. antiholomorphic) if \ $\alpha=J\alpha$ \ (resp. \ $\alpha$ \ is
perpendicular to \ $J\alpha$). A pair \ $\{X,Y\}$ \ of tangent vectors is said to be holomorphic
(resp. antiholomorphic) if \ $span \{X,Y\}$ \ is holomorphic (resp. antiholomorphic). 
Analogously are defined antiholo\-morphic triples of vectors. A 2-plane \ $\alpha$ \ in \ $T_pM$ \
or in \ $(T_pM)^\C$ \ is said to be weakly isotropic if the restriction of \ $g$ \ on \ $\alpha$ \
is of rank 1.

The curvature tensor \ $R$ \ of \ $M$ \ is defined by
$$
	R(X,Y)=[\nabla_X,\nabla_Y]-\nabla_{[X,Y]} \ ,
$$
where  $\nabla$  is the covariant differentiation on  $M$. Then the sectional curvature 
of a non\-degenerate 2-plane  $\alpha$, spaned by the orthonormal vectors \ $X,Y$ \ is 
defined as usual by
$$
	K(\alpha)=K(X,Y)=R(X,Y,Y,X)  \ .
$$

Let \ $R^\C$ \ and \ $\pi_1^\C$ \ be the natural extensions of \ $R$ \ and \ $\pi_1$, where
$$
	\pi_1(X,Y,Z,U)=g(X,U)g(Y,Z)-g(X,Z)g(Y,U) \ .
$$
Then we can consider also the sectional curvature
$$
	K^\C(\alpha) = R^\C(X,Y,Y,X)/\pi_1^\C(X,Y,Y,X) \ ,
$$
where \ $\alpha$ \ is a nondegenerate 2-plane spaned by the vectors \ $X,Y$ \ in \ $(T_pM)^\C$.

The manifold \ $M$ \ is said to be of pointwise constant holomorphic sectional curvature
if for any point \ $p$ \ in \ $M$ \ there exists a constant \ $c(p)$, such that for each unit
vector \ $X$ \ in \ $T_pM$ \ the holomorphic sectional curvature \ $H(X)=K(X,JX)$ \ does not depend
on \ $X$, i.e. \ $H(X)=c(p)$. Similarly are defined almost Hermitian manifolds of pointwise
constant antiholomorphic sectional curvature. Finally, \ $M$ \ is said to be of pointwise 
constant totally real biholomorphic sectional curvature if for any point \ $p\in M$ \
there exists a constant $c(p)$, such that for each antiholomorphic orthonormal pair $\{X,Y\}$
in \ $T_pM$ \ spaning a nondegenerate 2-plane, the totally real biholomorphic sectional curvature \
$H(X,Y)=R(X,JX,JY,Y)$ \ does not depend on \ $\{X,Y\}$, i.e. \ $H(X,Y)=c(p)$.

Let \ $M$ \ be a \ $2m$-dimensional almost Hermitian manifold of indefinite metric \ $g$ \ of signature \
$(2s,2(m-s)$, i.e. the tangent space is isometric to \ $\R_{2s}^{2m}$ \ with the inner product
$$
	<X,Y> = -\sum_{i=1}^{2s}X^iY^i + \sum_{j=2s+1}^{2m}X^jY^j  \ .
$$

A pair \ $\{ x,a \}$ \ of tangent vectors at a point \ $p\in M$ \ is said to be orthonormal of
signature $(+,-)$ if \ $g(x,x)=1$, $g(a,a)=-1$, $g(x,a)=0$. Orthonormal triples of certain
signature are determined in a similar way. Similarly, one speaks about the signature of a 2-plane 
\ $\alpha$ \ depending on the signature of the restriction of \ $g$ \ on \ $\alpha$ \ [3].

If \ $M$ \ is of definite metric we call a pair \ $\{ x,a \}$ \ of vectors in \ $(T_pM)^\C$ \ to be 
of signature $(+,-)$ if \ $g(x,x)=1$, $g(a,a)=-1$, $g(x,a)=0$. We note that there exist
nonisotropic vectors in \ $(T_pM)^\C$ \ which have no signature.

Along this paper the symbols \ $x,y,z$ \ will be used only for unit vectors of signature
$(+)$, the symbol \ $a$ - for unit vector of signature $(-)$, the symbols \ $X,Y,Z,U,V$ -
for unit vectors of arbitrary signature and the symbol \ $\xi$ - for isotropic vector.

\vspace{0.1 in}
2. \underline {Indefinite almost Hermitian manifolds}. We shall use the following lemma.

\underline{Lemma 1}. If \ $M$ \ is a $2m$-dimensional indefinite almost Hermitian manifold, 
$m>1$, satisfying
$$
	R(x,Jx,Jx,a)+R(x,Jx,Ja,x)=0   \leqno (1)
$$
for any antiholomorphic pair \ $\{x,a\}$, then \ $M$ is of pointwise constant holomorphic
sectional curvature.

\underline{Proof}. Let \ $\{x,a\}$ \ be an orthonormal antiholomorphic pair and \ $|t|<1$.
From (1) it follows 
$$
	R(x+ta,Jx+tJa,Jx+tJa,tx+a)+R(x+ta,Jx+tJa,tJx+Ja,x+ta)=0  \ .
$$
Since this equality holds for any $t$ with $|t|<1$, we obtain
$$
	2H(x)+2R(x,Jx,Ja,a)+2R(x,Ja,Jx,a)-K(x,Ja)-K(Jx,a)=0 \ ,  \leqno (2)
$$
$$
	2H(a)+2R(x,Jx,Ja,a)+2R(x,Ja,Jx,a)-K(x,Ja)-K(Jx,a)=0 \ ,  \leqno (3)
$$
$$
	R(a,Ja,Ja,x)+R(a,Ja,Jx,a)=0 \ .  \leqno (4)
$$
From (2) and (3) we find
$$
	H(x)=H(a) \ .  \leqno (5)
$$

If \ $m>2$ \ (5) implies that \ $M$ \ is of pointwise constant holomorphic sectional curvature.
If \ $m=2$ \ we can obtain the assertion using (1), (4) and (5).

\vspace{0.1in}
\underline{Theorem 1}. Let \ $M$ \ be a \ $2m$-dimensional almost Hermitian manifold of
indefinite metric of signature \ $(2s,2(m-s))$, $m>1$. If for each point \ $p\in M$ \
there exists a constant \ $c(p)$ \ such that for any vector \ $x \in T_pM$ \ the holomorphic
sectional curvature satisfies 
$$
	|H(x)| \le c(p) \ ,
$$ 
then \ $M$ \ is of pointwise constant holomorphic sectional curvature.

\underline{Proof}. Let \ $p\in M$ \ and \ $\{x,a\}$ \ be an arbitrary antiholomorphic 
orthonormal pair in \ $p$. Then from
$$
	| H( \frac{x+ta}{\sqrt{1-t^2}} ) | \le c(p)
$$
for \ $|t|<1$ \ we find 
$$
	\begin{array}{l}
		|H(x)+2t \{R(x,Jx,Jx,a)+R(x,Jx,Ja,x)\}  \\
		+t^2\{ 2R(x,Jx,Ja,a)+2R(x,Ja,Jx,a)-K(x,Ja)-K(Jx,a)\}  \\
		+2t^3\{ R(a,Ja,Ja,x)+R(a,Ja,Jx,a)\} +t^4H(a)| \le c(p)(1-t^2)^2 \ .
	\end{array} \leqno (6)
$$
Hence we obtain
$$
	\begin{array}{l}
		|H(x)		+t^2\{ 2R(x,Jx,Ja,a)+2R(x,Ja,Jx,a)\\
		-K(x,Ja)-K(Jx,a)\}   +t^4H(a)| \le c(p)(1-t^2)^2 \ .
	\end{array} \leqno (7)
$$
On the other hand, from (6) by continuity we get
$$
	\begin{array}{l}
		H(x)+2t \{R(x,Jx,Jx,a)+R(x,Jx,Ja,x)\}  \\
		+t^2\{ 2R(x,Jx,Ja,a)+2R(x,Ja,Jx,a)-K(x,Ja)-K(Jx,a)\}  \\
		+2t^3\{ R(a,Ja,Ja,x)+R(a,Ja,Jx,a)\} +t^4H(a) =0 
	\end{array} 
$$
for \ $t=\pm 1$, which implies
$$
	R(x,Jx,Jx,a)+R(x,Jx,Ja,x)+
		 		 R(a,Ja,Ja,x)+R(a,Ja,Jx,a)=0  \ . \leqno (8)
$$
From (6) and (8) we derive
$$
	\begin{array}{l}
		|H(x)+t^2\{ 2R(x,Jx,Ja,a)+2R(x,Ja,Jx,a)-K(x,Ja)-K(Jx,a)\}  \\
		+2t(1-t^2)\{ R(x,Jx,Jx,a)+R(x,Jx,Ja,x)\} +t^4H(a)| \le c(p)(1-t^2)^2 \ .
	\end{array} 
$$
Hence using (7) we get
$$
	|t\{R(x,Jx,Jx,a)+R(x,Jx,Ja,x)\}| \le c(p)(1-t^2) \ .
$$
By continuity we find (1) and the assertion follows from Lemma 1.

We shall use the following theorem.

\vspace{0.1in}
\underline{Theorem A}. [2] Let \ $M$ \ be a \ $2m$-dimensional almost Hermitian manifold of
indefinite metric of signature \ $(2s,2(m-s))$, $m>2$. If \ $R(X,\xi,\xi,X)=0$ \ whenever \
$span\{X,Y\}$ \ is a weakly isotropic antiholomorphic 2-plane, then \ $M$ \ is of pointwise
constant antiholomor\-phic sectional curvature. 

\underline{Proof}. Let e.g. \ $m-s>1$ \ and \ $\{x,y,a\}$ \ be an arbitrary orthonormal 
antiholomorphic triple in \ $T_pM$. From the condition we find \ 
$R(x,y+a,y+a,x)=0$ \ which implies
$$
	K(x,y)=K(x,a) \ ,  \leqno (9)
$$
$$
	R(x,y,a,x)=0  \ .  \leqno (10)
$$
Replacing in (9) \ $y$ \ by \ $(y+Jy)/\sqrt 2$, or \ $a$ \ by \ $(a+Ja)/\sqrt 2$, or \  $\{ x,y \}$ \ by \ 
$\{(x+y)/\sqrt 2$, $ (x-y)/\sqrt 2\}$ we get
$$
	R(x,y,Jy,x)=0\ , \quad R(x,a,Ja,x)=0 \ , \quad R(x,a,a,y)=0 \ ,  \leqno (11)
$$ 
respectively. Analogously, the first equality of (11) implies
$$
	R(a,y,Jy,a)=0  \ . \leqno (11')
$$

Let \  $X, Y, Z$ \ be arbitrary unit vectors in \ $T_pM$ \ so that 
\ $g(X,Y)=g(X,JY)=g(X,Z)=g(X,JZ)=0$ \ and \  $span\{X,Y\}$ \ and 
\ $span\{X,Z\}$ \ are nondegenerate. Applying (9),(10), (11) and (11') we obtain
$$
	K(X,Y)=K(X,Z) \ .   \leqno (12)
$$  

Suppose that \ $\alpha,\beta$ \ are both arbitrary nondegenerate antiholomorphic 
2-planes in \ $T_pM$ \ with orthonormal bases \ $\{X,Y\}$ \ and \ $\{Z,U\}$,
respectively. Let \ $V$ \ be a unit vector in \ $T_pM$ \ so that \ $V$ \ is normal to \ 
$span\{Y,JY,Z,JZ\}$. According to (12) we have \ $K(\alpha)=K(X,Y)=K(Y,V)=K(V,Z)=K(Z,U)=K(\beta)$.
Hence \ $M$ \ is of pointwise constant antiholomorphic sectional curvature.

\vspace{0.1in}
\underline{Theorem 2}. Let \ $M$ \ be a \ $2m$-dimensional almost Hermitian manifold of
indefinite metric of signature \ $(2s,2(m-s))$, $m>2$. If for each point \ $p\in M$ \ 
there exists a constant \ $c(p)$ \ such that the sectional curvature of any antiholomorphic
2-plane \ $\alpha$ \ in \ $T_pM$ \ is bounded from above by \ $c(p)$, i.e.
$$
	K(\alpha) \le c(p) \ ,
$$
then \ $M$ \ is of pointwise constant antiholomorphic sectional curvature.

\underline{Proof}. Let e.g. \ $m-s \ge 2$. We choose an orthonormal antiholomorphic
triple \ $\{x,y,a\}$ \ in \ $T_pM$. Then the inequality \ 
$K(span\{x+ta,y\}) \le c(p)$ \ for \ $t \ne \pm 1$ \ implies
$$
	R(x+ta,y,y,x+ta) \
	\left\{  
		\begin{array}{lcr}
			\le c(p)(1-t^2)   & {\rm if} & |t|<1 \ ,  \\
			\ge c(p)(1-t^2)   & {\rm if} & |t|>1 \ .
		\end{array}
	\right.
$$
Hence by continuity we find \ $R(x+a,y,y,x+a)=0$ \ and using Theorem A we obtain
the assertion.

\underline{Remark 1}. The requirement \ $K(\alpha)\le c(p)$ \ in Theorem 2 can be 
replaced by \ $K(\alpha)\ge c(p)$. The conclusion rest true also if \ 
$|K(\alpha)|\le c(p)$ \ for any antiholomorphic 2-plane \ $\alpha$ \ of signature \ 
$(+,-)$. Analogously for antiholomorphic 2-planes of signature (+,+) or of
signature $(-,-)$.

\vspace{0.1in}
\underline{Theorem 3}. Let \ $M$ \ be a \ $2m$-dimensional almost Hermitian manifold of
indefinite metric of signature \ $(2s,2(m-s))$, $m>2$. If \ $R(X,JX,J\xi,\xi)=0$ \ 
whenever \ $span\{X,\xi\}$ \ is a weakly isotropic antiholomorphic 2-plane, 
then \ $M$ \ is of pointwise constant totally real biholomorphic sectional curvature.

\vspace{0.1in}
\underline{Theorem 4}. Let \ $M$ \ be a \ $2m$-dimensional almost Hermitian manifold of
indefinite metric of signature \ $(2s,2(m-s))$, $m>2$. If for each point \ $p\in M$ \ 
there exists a constant \ $c(p)$ \ such that for any antiholomorphic 2-plane \ $\alpha$ \ 
in \ $T_pM$ \ the totally real biholomorphic sectional curvature is bounded from above,
then \ $M$ \ is of pointwise constant totally real biholomorphic sectional curvature.

Theorem 3 and Theorem 4 can be proved in a similar way as Theorem A and Theorem 2, so
we omit the proofs. A remark, like Remark 1 holds for Theorem 4 too.

\vspace{0.2in}
3. \underline{Definite almost Hermitian manifolds}.

\vspace{0.1in}
\underline{Lemma 2}. If \ $M$ \ is a definite \ $2m$-dimensional  almost Hermitian manifold, 
$m>1$, satisfying
$$
	R(x,Jx,Jx,y)+R(x,Jx,Jy,x)=0   
$$
for any antiholomorphic pair \ $\{x,y\}$, then \ $M$ \ is of pointwise constant holomorphic
sectional curvature.

The proof is essentially the same as in the indefinite case.

\vspace{0.1in}
\underline{Theorem 5}. Let \ $M$ \ be a \ $2m$-dimensional almost Hermitian manifold of
definite metric,  $m>1$. If for each point \ $p\in M$ \ 
there exists a constant \ $c(p)$ \ such that the absolute value of the real or of the 
imaginary part of the holomorphic sectional curvature \ $H^\C(\alpha)$ \ of each 
holomorphic 2-plane \ $\alpha$ \ in \ $(T_pM)^\C$ \ of signature (+,+) is bounded by \ 
$c(p)$, then \ $M$ \ is of pointwise constant holomorphic sectional curvature.

\underline{Proof}. Suppose that e.g. \ $|Re\, H^\C(\alpha)|$ \ is bounded. Let \ 
$\{x,y\}$ \ be an arbitrary orthonormal antiholomorphic pair in \ $T_pM$. Then from
$$
	|Re\, H^\C\left(\frac{x+ity}{\sqrt{1-t^2}}\right)| \le c(p)
$$
for \ $|t|<1$, we obtain
$$
	\begin{array}{l}
		|H(x)-t^2\{K(x,Jy)+2R(x,Jx,Jy,y)+2R(x,Jy,Jx,y)\\
		+K(Jx,y)\}+t^4H(y)| \le c(p)(1-t^2)^2 \ .  
	\end{array}\leqno (12)
$$
Hence by continuity we have
$$
	H(x)+H(y) = K(x,Jy)+2R(x,Jx,Jy,y)+2R(x,Jy,Jx,y)	+K(Jx,y) \ .  \leqno (13)
$$
From (12) and (13) it follows \ $|H(x)-t^2H(y)| \le c(p)(1-t^2) $ \ for
all \ $t$, $|t|<1$. Now by continuity we find
$$
	H(x)=H(y) \ .   \leqno (14)
$$ 

If \ $m>2$, then \ $M$ \ is of pointwise constant holomorphic sectional curvature
because of (14). Let \ $m=2$. From (13) and (14) we obtain
$$
	2H(x)= K(x,Jy)+2R(x,Jx,Jy,y)+2R(x,Jy,Jx,y)	+K(Jx,y) \ .  \leqno (15)
$$
Replacing \ $\{x,y\}$ \ by  \ $\{(x+ty)/\sqrt{1-t^2},\ (tx-y)/\sqrt{1-t^2}\}$, \ $|t|<1$ \  
in (15) we get
$$
	5R(x,Jx,Jx,y)+5R(x,Jx,Jy,x)=3R(x,Jy,Jy,y)+3R(y,Jy,Jx,y) \ ,
$$
$$
	3R(x,Jx,Jx,y)+3R(x,Jx,Jy,x)=5R(x,Jy,Jy,y)+5R(y,Jy,Jx,y) \ .
$$
These equalities give \ $ R(x,Jx,Jx,y)+R(x,Jx,Jy,x)=0$ \ and the assertion
follows from Lemma 2. 

\vspace{0.1in}
\underline{Lemma 3}. If \ $M$ \ is a  \ $2m$-dimensional definite almost Hermitian manifold, 
$m>2$. Then the following conditions are equivalent:

a) $R(x,y,z,x)=0$ \ for an arbitrary orthonormal antiholomorphic triple \ $\{x,y,z\}$;

b) $K(x,y)=K(x,z)$ \ for an arbitrary orthonormal antiholomorphic triple \ $\{x,y,z\}$;

c) $M$ \ is of pointwise constant antiholomorphic
sectional curvature.

\underline{Proof}. We shall show  that a) and b) are equivalent. In fact, if we replace
\ $\{y,z\}$ \ by \ $ \{(y-z)/\sqrt2,(y+z)/\sqrt2\}$ \ in a) we obtain b). Similarly b)
implies a). The implication a) $\rightarrow$ c) (or b) $\rightarrow$ c)) can be obtained 
by using the same arguments as in the proof of Theorem A. Finally, c) obviously implies a). 

\vspace{0.1in}
\underline{Theorem 6}. Let \ $M$ \ be a \ $2m$-dimensional almost Hermitian manifold of
definite metric,  $m>2$. If \ $R^\C(x,\xi,\xi,x)=0$, whenever \ $span\{x,\xi\}$ \ is 
a weakly isotropic antiholomorphic 2-plane, then \ $M$ \ is of pointwise constant 
antiholomorphic sectional curvature.

\underline{Proof}. Let \ $\{x,y,z\}$ \ be an arbitrary orthonormal antiholomorphic triple.
Then \ $R^\C(x,y+iz,y+iz,x)=0$ \ implies \ $K(x,y)=K(x,z)$ \ and the assertion follows from Lemma 3.

\vspace{0.1in}
\underline{Theorem 7}. Let \ $M$ \ be a \ $2m$-dimensional almost Hermitian manifold of
definite metric,  $m>2$. If for each point \ $p\in M$ \ 
there exists a constant \ $c(p)$ \ such that for each antiholomorphic nondegenerate 
2-plane \ $\alpha$ \ the real or the imaginery part of \ $K^\C(\alpha)$ \  is bounded 
either form above or from below by \ $c(p)$, then \ $M$ \ is of pointwise constant 
antiholomorphic sectional curvature.

\underline{Proof}. Suppose that e.g. the boundedness is from below. Let \ $\{x,y,z\}$ \   
be an arbitrary orthonormal antiholomorphic triple in \ $T_pM$. Then from
$$
	Re\, K^\C\left(x,\frac{y+itz}{\sqrt{1-t^2}}\right) \ge c(p)
$$
for real \ $t\ne\pm 1$, we obtain
$$
	K(x,y)-t^2K(x,z) \
	\left\{  
		\begin{array}{lcr}
			\ge c(p)(1-t^2)   & {\rm if} & |t|<1 \ ,  \\
			\le c(p)(1-t^2)   & {\rm if} & |t|>1 \ .
		\end{array}
	\right.
$$
By continuity we get the condition b) of Lemma 3 and so the assertion follows from Lemma 3.

Because of the analogy with the case of totally real biholomorphic sectional
curvature we shall not formulate the corresponding results.

\vspace{0.7in}
\centerline{\large REFERENCES}

\vspace{0.1in}
\noindent
1. M. Barros, A. Romero. Indefinite Kaehlerian Manifolds. 
Math. Ann., \underline{261}(1982), 

55-62.

\noindent
2. A. Borisov, G. Ganchev, O. Kassabov. Curvature Properties and Isotropic Planes
of 

Riemannian and Almost Hermitian manifolds of Indefinite Metrics.
Ann. Univ. Sof., 

Fac. Math.  M\'ec., \underline{78}(1984), 53-56.

\noindent
3. M. Dajczer, K. Nomizu. On Sectional Curvature of Indefinite Metrics.
Math. Ann., 

\underline{247}(1980), 279-282.

\noindent
4. G. Ganchev, V. Mihova. On the Conformal Curvature Tensor in the Riemannian and 

the Almost Hermitian Geometry.  Ann. Univ. Sofia, Math., \underline{78}(1984), 158-171.

\noindent
5. O. Kassabov, A. Borisov. On the Sectional Curvature of K\"ahler Manifolds of Indefinite 

Metrics. PLISKA Studia mathematica bulgarica, \underline9(1987), 66-70.

\noindent
6. R. S. Kulkarni. The Values of Sectional Curvature of Indefinite Metrics.
Comment. 

Math. Helvetici, \underline{54}(1979), 173-176.

\end{document}